Stanislav Drastich, mailto:mydrass@hotmail.com

# Rapid growth sequences

## *Introduction*

We are going to explain some basic rapid growth sequences having important properties. Studying Fermat sequence

$$2^{2^n} + 1$$

we can simply find important property of this type of sequences that might be described as follows

$$\left(2^{2^n} + 1\right) = \left(2^{2^0} + 1\right)\left(2^{2^1} + 1\right)...\left(2^{2^{n-1}} + 1\right) + 2$$

This is some type of rapid growth sequence having more important common properties based on sequence construction.
Thus we are on duty to ask question if there are some other similar sequences having the same properties.
This approach leads to the arranged ideas we are going to introduce now.

## *Rapid growth sequences construction*

In fact, we are going to explain common method of construction of rapid growth sequences.
For this reason, let us choose any natural number $\nabla$ and let us fix this number.

### Definition 1

Let $\nabla$ be any fix natural number.
Let $\Delta_0 = 1 + \nabla$
$$\Delta_n = \Delta_{n-1}^2 - \nabla \Delta_{n-1} + \nabla$$
Then we call the sequence $\Delta_0, \Delta_1, \Delta_n, ...$ rapid growth sequence of first type with $\nabla$ base.

### Definition 2

Let $\nabla$ be any fix natural number.
Let $\Delta_0 = 1 - \nabla$
$$\Delta_n = \Delta_{n-1}^2 + \nabla \Delta_{n-1} - \nabla$$
Then we call the sequence $\Delta_0, \Delta_1, \Delta_n, ...$ rapid growth sequence of second type with $\nabla$ base.

Thus we have constructed infinitely many different sequences of properties similar to the properties of Fermat sequence, as we are going to show. Let us summarize general valuable properties.

Stanislav Drastich, mailto:mydrass@hotmail.com

## *Basic characterization of rapid growth sequences*

### Lemma 1

Let $\Delta_0, \Delta_1, \Delta_n, \ldots$ be rapid growth sequence of first (second) type.
Then $\Delta_n = \Delta_0 \Delta_1 \ldots \Delta_{n-2} \Delta_{n-1} \pm \nabla$

### *Proof*

Because of $\Delta_0, \Delta_1, \Delta_n, \ldots$ is rapid growth sequence of first (second) type, thus we have

$\Delta_1 = \Delta_0^2 \mp \nabla \Delta_0 \pm \nabla$

$\Delta_1 = (1 \pm \nabla)^2 \mp \nabla(1 \pm \nabla) \pm \nabla$

$\Delta_1 = 1 \pm 2\nabla + \nabla^2 - \nabla^2 \mp \nabla \pm \nabla$

$\Delta_1 = 1 \pm 2\nabla$

But $\Delta_0 \pm \nabla = (1 \pm \nabla) \pm \nabla = 1 \pm 2\nabla = \Delta_1$

Let $\Delta_n = \Delta_0 \Delta_1 \ldots \Delta_{n-2} \Delta_{n-1} \pm \nabla$.

Then under the definition 1 (2 respectively) $\Delta_{n+1} = \Delta_n^2 \mp \nabla \Delta_n \pm \nabla$ partially exchanging appropriate indexes. Thus

$\Delta_{n+1} = \Delta_n^2 \mp \nabla \Delta_n \pm \nabla$

$\Delta_{n+1} \mp \nabla = \Delta_n^2 \mp \nabla \Delta_n$

$\Delta_{n+1} \mp \nabla = (\Delta_n \mp \nabla)\Delta_n$

Under the assumption $\Delta_n = \Delta_0 \Delta_1 \ldots \Delta_{n-2} \Delta_{n-1} \pm \nabla$ we have

$\Delta_{n+1} \mp \nabla = (\Delta_0 \Delta_1 \ldots \Delta_{n-2} \Delta_{n-1} \pm \nabla \mp \nabla)\Delta_n$

$\Delta_{n+1} \mp \nabla = \Delta_0 \Delta_1 \ldots \Delta_{n-2} \Delta_{n-1} \Delta_n$

$\Delta_{n+1} = \Delta_0 \Delta_1 \ldots \Delta_{n-2} \Delta_{n-1} \Delta_n \pm \nabla$

q.e.d.

### Examples

| $\nabla$ | Rapid growth sequence of first type | Rapid growth sequence of second type |
|---|---|---|
| 1 | 2, 3, 7, 43, 1807, 3263443, ... | 0, -1, -1, -1, -1, -1, ... |
| 2 | 3, 5, 17, 257, 65537, 4294967297, ... | -1, -3, 1, 1, 1, 1, ... |
| 3 | 4, 7, 31, 871, 756031, 571580604871, ... | -2, -5, 7, 67, 4687, 21982027, ... |
| 4 | 5, 9, 49, 2209, 4870849, 23725150497409, ... | -3, -7, 17, 353, 126017, 15880788353, ... |
| 5 | 6, 11, 71, 4691, 21982031, 483209576974811, ... | -4, -9, 31, 1111, 1239871, 1537286295991, ... |
| 6 | 7, 13, 97, 8833, 77968897, 6079148431583233, ... | -5, -11, 49, 2689, 7246849, 52516863909889, ... |
| 7 | 8, 15, 127, 15247, 232364287, | -6, -13, 71, 5531, 30630671, |

Stanislav Drastich, mailto:mydrass@hotmail.com

|    |                                                              |                                                        |
|----|--------------------------------------------------------------|--------------------------------------------------------|
|    | 53993160246468367, ...                                       | 938238220324931, ...                                   |
| 8  | 9, 17, 161, 24641, 606981761, 368426853330807041, ...        | -7, -15, 97, 10177, 103652737, 10743890716813057, ...  |
| 9  | 10, 19, 199, 37819, 1429936399, 2044718092315659619, ...     | -8, -17, 127, 17263, 298166527, 88903280506740463, ... |
| 10 | 11, 21, 241, 55681, 3099816961, 9608865160705105921, ...     | -9, -19, 161, 27521, 757680641, 5740799613229772 81, ...|

For instance, choosing $\nabla = 2$ for rapid growth sequence of first type we obtain Fermat sequence as introduced by $2^{2^{n-1}} + 1$ formula, where n is natural number.

Rapid growth sequences introduced in this manner have important properties related to prime divisibility characterization we are going shortly to describe.

## Lemma 2

Let $\Delta_0, \Delta_1, \Delta_n, ...$ be a rapid growth sequence of first (second) type with $\nabla$ base. Then common prime divisor of every two different $\Delta_n, \Delta_{n+m}$ where $1 \leq m$, is also prime divisor of $\nabla$.

### *Proof*

Let us have two different members of rapid growth sequence of first (second) type. Thus we can suppose to have chosen $\Delta_n, \Delta_{n+m}$, where $1 \leq m$. In respect to the statement of Lemma 1 we can write down $\Delta_{n+m} = \Delta_0 \Delta_1 ... \Delta_{n-1} \Delta_n \Delta_{n+1} ... \Delta_{n+m-1} \pm \nabla$. Let p be their common prime divisor, in other words $p \mid \Delta_n$ and $p \mid \Delta_{n+m}$. But we have supposed, that $1 \leq m$, so there is $\Delta_n$ involved in $\Delta_0 \Delta_1 ... \Delta_{n-1} \Delta_n \Delta_{n+1} ... \Delta_{n+m-1}$, that means the left side of equation $\Delta_{n+m} = \Delta_0 \Delta_1 ... \Delta_{n-1} \Delta_n \Delta_{n+1} ... \Delta_{n+m-1} \pm \nabla$ is divisible by p and at the right side of the same equation we have sum of number divisible by p and $\nabla$. Thus $p \mid \nabla$.
q.e.d.

## Lemma 3

Let $\Delta_0, \Delta_1, \Delta_n, ...$ be a rapid growth sequence of first (second) type with $\nabla$ base. Then for any prime divisor $p_\nabla$ of $\nabla$ follows
$\Delta_n \equiv 1 (\mod p_\nabla)$

### *Proof*

Let us have a rapid growth sequence of first (second) type with $\nabla$ base.
Then $\Delta_0 = 1 \pm \nabla \equiv 1 (\mod p_\nabla)$.
Thus let $\Delta_n \equiv 1 (\mod p_\nabla)$.
Following sequence member is $\Delta_{n+1} = \Delta_n^2 \mp \nabla \Delta_n \pm \nabla$

Stanislav Drastich, mailto:mydrass@hotmail.com

$$\Delta_{n+1} \equiv \Delta_n^2 \mp \nabla \Delta_n \pm \nabla \equiv 1^2 \mp \nabla \pm \nabla \equiv 1 (\bmod p_\nabla)$$

Thus $\Delta_{n+1} \equiv 1 (\bmod p_\nabla)$.

q.e.d.

## Theorem 1

Let $\Delta_0, \Delta_1, \Delta_n, \ldots$ be a rapid growth sequence of first (second) type with $\nabla$ base.
Then every two different $\Delta_n, \Delta_m$ where $n \neq m$ are relatively prime.

### *Proof*

Let $\Delta_0, \Delta_1, \Delta_n, \ldots$ be a rapid growth sequence of first (second) type with $\nabla$ base.
Let us suppose, that there exists some prime divisor p of two different sequence numbers $\Delta_n, \Delta_m$, where $n \neq m$. In respect to the statement of Lemma 2 p is also prime divisor of $\nabla$ as well.
On the other hand following result of Lemma 3, for every $\Delta_n \equiv 1 (\bmod p)$ that is in controversy to the assumption $\Delta_n \equiv 0 (\bmod p)$ and $\Delta_m \equiv 0 (\bmod p)$.

q.e.d.

## Corollary 1

Let $\Delta_0, \Delta_1, \Delta_n, \ldots$ be a rapid growth sequence of first (second) type with $\nabla$ base.
Let $p_\nabla$ be any prime divisor of $\nabla$. Let n be zero or natural number, let m be natural number.
Then $p_\nabla \mid (\Delta_n \ldots \Delta_{n+m-1} - 1)$.

### *Proof*

Let us choose two different numbers $\Delta_n, \Delta_{n+m-1}$. In respect to Lemma 3 $\Delta_n \equiv 1 (\bmod p_\nabla)$ and $\Delta_{n+m-1} \equiv 1 (\bmod p_\nabla)$ are of the same residue classes modulo $p_\nabla$. Thus $p_\nabla$ is prime divisor of their difference $\Delta_m - \Delta_n$.

$$\Delta_{n+m-1} - \Delta_n = (\Delta_0 \ldots \Delta_{n-1} \Delta_n \Delta_{n+1} \ldots \Delta_{n+m-1} \pm \nabla) - (\Delta_0 \ldots \Delta_{n-1} \pm \nabla)$$

$$\Delta_{n+m-1} - \Delta_n = (\Delta_0 \ldots \Delta_{n-1} \Delta_n \Delta_{n+1} \ldots \Delta_{n+m-1} - \Delta_0 \ldots \Delta_{n-1}) \pm \nabla \mp \nabla$$

$$\Delta_{n+m-1} - \Delta_n = \Delta_0 \ldots \Delta_{n-1} (\Delta_n \Delta_{n+1} \ldots \Delta_{n+m-1} - 1)$$

The left side of mentioned equation is divisible by $p_\nabla$, thus the right side has to be divisible by $p_\nabla$ as well. We have chosen $p_\nabla$ as a prime number, so at least one of the factors $\Delta_0, \ldots, \Delta_{n-1}, or (\Delta_n \Delta_{n+1} \ldots \Delta_{n+m-1} - 1)$ has to be divisible by $p_\nabla$. In respect to the statement of Lemma 3 for every natural n and every prime divisor $p_\nabla \mid \nabla$ is $\Delta_n \equiv 1 (\bmod p_\nabla)$. That means, that we have just one possible factor divisible by $p_\nabla$ and it is $(\Delta_n \ldots \Delta_{n+m-1} - 1)$.

Stanislav Drastich, mailto:mydrass@hotmail.com

Thus $p \mid (\Delta_i ... \Delta_{j-1} - 1)$ should be true.
q.e.d.

In other words, we can say that difference of any two numbers of rapid growth sequence is divisible by prime divisors of $\nabla$.

### *Example*
For $\nabla = 5$
71*4691*21982031-1=7321357226890=5*1464271445378

For $\nabla = 10$
27521*757680641-1 = 20852128920960 = 2*10426064460480 = 5*4170425784192

In fact, we can say much more strongly, that any product of numbers of rapid growth sequence is "nearly divisible" by prime divisors of $\nabla$ as we will show now.

## Corollary 2
Let $\Delta_0, \Delta_1, \Delta_n, ...$ be a rapid growth sequence of first (second) type with $\nabla$ base.
Let $p_\nabla$ be any prime divisor of $\nabla$.
Let $i_1, i_2, ..., i_m$ be natural numbers or zeroes (indexes).
Then $p_\nabla \mid (\Delta_{i_1} ... \Delta_{i_m} - 1)$.

### *Proof*
We have chosen m (not necessarily different) indexes to rapid growth sequence of first (second) type with $\nabla$ base and prime divisor $p_\nabla$ of $\nabla$.
According to the result of Lemma 3 for any index and for any prime divisor $p_\nabla$ of $\nabla$ is
$\Delta_n \equiv 1 (\mod p_\nabla)$.
Thus $\Delta_{i_1} ... \Delta_{i_m} \equiv 1.1...1 \equiv 1 (\mod p_\nabla)$
$\Delta_{i_1} ... \Delta_{i_m} - 1 \equiv 0 (\mod p_\nabla)$
In other words, $p_\nabla \mid (\Delta_{i_1} ... \Delta_{i_m} - 1)$
q.e.d.

## Corollary 3
Fermat sequence numbers are relatively prime.

### *Proof*
Let us choose $\nabla = 2$. Let us construct rapid growth sequence of first type with 2 base.
Then

Stanislav Drastich, mailto:mydrass@hotmail.com

$\Delta_0 = 1 + 2 = 3 = 2^{2^0} + 1$

Let $\Delta_n = 2^{2^n} + 1$. Then $\Delta_{n+1} = \Delta_n^2 - 2\Delta_n + 2$

$\Delta_{n+1} = \left(2^{2^n} + 1\right)^2 - 2\left(2^{2^n} + 1\right) + 2$

$\Delta_{n+1} = 2^{2^n} \cdot 2^{2^n} + 2 \cdot 2^{2^n} + 1 - 2 \cdot 2^{2^n} - 2 + 2$

$\Delta_{n+1} = 2^{2^n} \cdot 2^{2^n} + 1$

$\Delta_{n+1} = 2^{2^n + 2^n} + 1$

$\Delta_{n+1} = 2^{2 \cdot 2^n} + 1$

$\Delta_{n+1} = 2^{2^{n+1}} + 1$

We have shown, that this sequence is on one side exactly the same sequence as Fermat sequence, and on the other side this sequence has been brought to life as a rapid growth sequence of first type with base 2 base.

Result of this corollary then follows directly from Theorem 1.

q.e.d.